\documentclass[%
sn-basic]{sn-jnl}% referee option is meant for double line spacing

\usepackage{graphicx}%
\usepackage{multirow}%
\usepackage{amsmath,amssymb,amsfonts}%
\usepackage{amsthm}%
\usepackage{mathrsfs}%
\usepackage[title]{appendix}%
\usepackage{xcolor}%
\usepackage{textcomp}%
\usepackage{manyfoot}%
\usepackage{booktabs}%
\usepackage{algorithm}%
\usepackage{algorithmicx}%
\usepackage{algpseudocode}%
\usepackage{listings}%
\usepackage{xspace}%

\newtheorem{example}{Example}

\newcommand{\PC}{\mbox{{\textrm{PC}}}\xspace}
\newcommand{\pc}{\PC}

\newcommand{\pr}{\Pr}
\newcommand{\bM}{\mbox{\bf M}}
\newcommand{\bR}{\mbox{\bf R}}

\newcommand{\tabref}[1]{\mbox{Table~\ref{tab:#1}}}
\newcommand{\figref}[1]{\mbox{Figure~\ref{fig:#1}}}

\newcommand{\secref}[1]{\mbox{\S$\,$\ref{sec:#1}}}
\newcommand{\Secref}[1]{\mbox{Section~\ref{sec:#1}}}
\renewcommand{\eqref}[1]{\mbox{(\ref{eq:#1})}}
\newcommand{\exref}[1]{\mbox{Example~\ref{ex:#1}}}

\newcommand{\cip}{\mbox{$\,\perp\!\!\!\perp\,$}}
\newcommand{\rr}{\mbox{\textup{\textrm{RR}}}\xspace}

\newcommand{\eg}{{\em e.g.\/}\xspace}

\begin{document}

\title[Individual Causation with Biased Data]{Individual Causation with Biased Data}


\author*[1]{\fnm{Monica} \sur{Musio}}\email{mmusio@unica.it}

\author[2]{\fnm{Philip} \sur{Dawid}}\email{apd@statslab.cam.ac.uk}
\equalcont{These authors contributed equally to this work.}


\affil[1]{\orgdiv{Dipartimento di Matematica ed Informatica}, \orgname{Universit\`a degli studi di Cagliari}, \orgaddress{\street{Via Ospedale, 72}, \city{Cagliari}, \postcode{09124}, \state{Sardinia}, \country{Italy}}}

\affil[2]{\orgdiv{Statistical Laboratory}, \orgname{University of
    Cambridge},

   \orgaddress{
     \street{
     %  Wilberforce Road
     }%,
     \city{
       % Cambridge
     }%,
     \postcode{
       % CB3 0WB},
       }
% \state{State},

 \country{UK}}}


\abstract{We consider the problem of assessing whether, in an
  individual case, there is a causal relationship between an observed
  exposure and a response variable.  When data are available on
  similar individuals we may be able to estimate prospective
  probabilities, but even under ideal consitions these are typically
  inadequate to identify the ``probability of causation'': instead we
  can only derive bounds for this.  These bounds can be improved or
  amended when we have information on additional variables, such as
  mediators or covariates.  When a covariate is unobserved or ignored,
  this will typically lead to biased inferences.  We show by examples
  how serious such biases can be.}




\maketitle

\section{Introduction}\label{sec:intro}

The field of causal inference subsumes two quite distinct general
activities: inference about ``Effects of Causes" (EoC), and inference
about ``Causes of Effects" (CoE) \citep{apd/mm:eoccoe}.  These require
quite different mathematical frameworks and analyses \citep{apd:kent}.

To illustrate the difference, consider the relationship between being
vaccinated against COVID-19, and developing myocarditis (within 28
days of vaccination).  We examine first the case of 20-year-old Sven,
who is considering getting the mRNA-1273 vaccine against COVID-19, but
is concerned about the possibility that he will develop myocarditis.
The relevant information for Sven is contained in the probabilities
that he will develop myocarditis, if he does, or if he does not, take
the vaccine.  In particular, the estimated probability of myocarditis
if vaccinated is (according to the data analysis of \citet{karlstad})
approximately $3 \times 10^{-4}$.  Here we have an EoC question and
answer: What might be the effect (possible myocarditis), consequent on
an applied cause (vaccination)?  EoC focuses on the predictive power
of causation, helping us anticipate outcomes and make informed
decisions regarding individual health, safety, and well-being.  And it
can be addressed by applying the standard tools of statistical
inference and decision theory, using relevant data, even taking into
account certain possible biases in the data, such as confounding
\citep{apd:found,pearl:book}.

Now our focus shifts to Anders, another 20-year-old, who did receive
the mRNA-1273 vaccine and unfortunately developed myocarditis.  The
question that arises here is not about predicting the outcome, which
is known, but rather about attributing causality: was it the vaccine
that caused the myocarditis?  This is what we mean by a CoE question.
And to address it we need to expand the mathematical framework, as
will be described below.  On doing so we find that (again based on the
data of \citet{karlstad}) the probability that Anders' myocarditis
was, indeed, a direct consequence of his vaccination (the
``probability of causation'', \pc) is estimated to be at least $0.97$.
Such a bound on \pc is typical in CoE problems, where, even with ideal
data, it is typically not possible to identify \pc exactly.  It is
striking how different the answers to EoC and CoE questions can be.

In this paper we concentrate attention on the CoE problem, where we
seek to understand whether there is a causal link between a putative
cause and an observed effect in an individual case.  We pay particular
attention to the potential for bias, and its consequences.  The data
at hand may be derived from an observational study, which inherently
carries the risk of various biases such as selection bias, confounding
variables, and incomplete information.  In this work we show, in
simple cases, how the presence of confounding bias in the data can
lead to completely misleading conclusions about the probability of
causation.

The paper is organised as follows.  In \secref{coe} we frame the CoE
problem as a counterfactual query, and in \secref{po} show how this
may be addressed using the device of potential outcomes.
\Secref{data} discusses just what can be ideally estimated from data,
which necessarily falls short of what is desired.  \Secref{basic}
explains how, in a simple case, estimable probabilities can be used to
set interval bounds on \pc.  \Secref{med} introduces an additional
variable that mediates the causal relationship: this can lead to
refined bounds on \pc.  By contrast, in \secref{suffcov} we introduce
a covariate, that can affect both exposure and response.  When it is
observed in the data, though not necessarily for Anders, we obtain new
bounds for \pc; but if it is not observed, and ignored, it acts as a
confounder, leading to a biased analysis.  By a simple example we show
how misleading this can be.  In \secref{medsuff} we combine the
previous two cases: again, including the mediator can refine the
bounds on \pc, but ignoring the covariate can be highly misleading.
Finally in \secref{conc} we review our results.

%\section{Vaccination and myocarditis}\label{sec:vacc}
%An analysis of data on
 % the populations of 4 Nordic countries (Karlstad et al., 2022)
 % estimated the probability of an individual vaccinated against
  %SARS-CoV-2 developing myocarditis within 28 days, classified by
  %vaccine, age, sex, and several other covariates.

   %For example, for a male
 % aged 16--24 receiving mRNA-1273, the rate was 3.7 per thousand person-years of
%  follow-up, compared with 0.19 if unvaccinated---an incidence rate
%  ratio (IRR) of 35.

%The IRRs were attenuated after adjusting for other covariates.

%\section{Two perspectives}\label{sec:persp}
%\begin{description}

%\item[Effects of Causes (EoC):]  

%Sven, aged 20, is about to
%  have the mRNA-1273 vaccine against COVID-19.    What is the
  %probability  that he will develop myocarditis (within 28 days)?
%  $$3 \times 10^{-4}$$
  
%\item [Causes of Effects (CoE):] Anders, aged 20, had the  mRNA-1273
%  vaccine, and developed myocarditis.  What is the probability that
%  this was {\em caused\/} by the  vaccine?

 % {Equivalently:  Given the above information, what is the probability he would not have developed
 %   myocarditis if he had  not had the vaccine?}
%  $$\geq \, 0.97$$
  
%\end{description}

\section{Causes of effects}\label{sec:coe}

Consider the scenario in which Anders was vaccinated and developed
myocarditis.  We would like to know whether it was the vaccination
that caused the myocarditis.  It will not be possible to decide this
fact for certain, but the aim is to estimate its probability, on the
basis of available evidence.

This is the kind of case that arises in ``toxic tort'' law suits,
where an individual claims damages from a company on the basis that it
was their drug that caused his adverse reaction.  How are we to
understand and investigate such a claim?  The usual way in which a
court frames this question is by the ``but for'' criterion: the
relationship between exposure and outcome (both known to have
occurred) is considered {\em causal\/} just when the outcome would not
have occurred, {\em but for\/} the exposure.  That is to say, if the
exposure had not been present, then the outcome would not have
happened.  For, if the outcome would have happened even in the absence
of the exposure, then it can not be attributed to the exposure.

Note that this approach requires consideration of a {\em
  counterfactual\/} scenario: we need to think about events in a
parallel universe where the known fact of exposure is negated.  This
leads us into murky philosophical waters.  In particular, since
neither counterfactual events nor the ``fact'' of causation can ever
be observed, somewhat subtle statistical analysis is required to
define the ``probability of causation'', \pc, in general, and to
estimate Anders's own individual probability of causation, $\pc_A$, in
particular, using data from similar individuals.
 
%\begin{center}
%{Binary exposure} $E$  \quad {Binary response} $R$ \\
%\end{center}
%\begin{itemize}
%\item Anders has been subjected to exposure (e.g. vaccination)  $E$: $E_A=1$

 % \item Anders develops undesirable outcome (e.g. myocarditis)  $R$: $R_A=1$ 
  
%  \item Interest in whether it was Anders's exposure that {caused}
 %   his outcome
%  \end{itemize}
   
%\subsubsection{Causes of effects: the Court perspective}\label{sec:court}

\section{Potential outcomes}\label{sec:po}

Statisticians have developed an approach based on ``potential
outcomes'' \citep{dbr:jep,dbr:as} to define \pc
\citep{tian/pearl:probcaus,apd/mm/rm}.  For a generic individual, let
binary variables $E$ and $R$ denote, respectively, exposure and
outcome (response) status.  We duplicate $R$ by introducing
``potential outcomes'':
\begin{itemize}
\item $R(0)$ denotes the outcome that would realise if $E$ takes
  value $0$ (no exposure)
\item $R(1)$ denotes the outcome that would realise if $E$ takes
  value $1$ (exposure)
\end{itemize}
(It is assumed that the potential outcome does not depend on how the
individual is exposed, whether by an external intervention or in an
undisturbed natural way.)
  
The actual outcome if exposed, $E=1$, is $R=R(1)$, which can be
observed.  In this case $R(0)$, the response if unexposed, becomes a
counterfactual variable, and can not be observed.  The impossibility
of observing both $R(0)$ and $R(1)$ in the same individual case has
been termed ``the fundamental problem of causal inference''
\citep{pwh:jasa}.
    
In spite of this problem, the conception is that both $R(0)$ and
$R(1)$ exist, simultaneously, even prior to realisation of $E$.  And
$(E, R(0), R(1))$ are considered to have a joint probability
distribution.  However, it is never possible to obtain samples from
this distribution.  At best, we can observe $(E,R) = (E,R(E))$.

\subsection{Counterfactual causation and \pc}
\label{sec:count-caus}

In terms of potential responses, we interpret ``exposure causes
outcome'' as the event $C=(R(0) = 0, R(1)=1)$: that is, the outcome
occurs if and only if the individual is exposed.

Knowing that Anders was vaccinated ($E_A=1$) and developed myocarditis
($R_A=1$), his probability of causation is thus
\begin{eqnarray}
  \pc_A &=& \pr(C_A \mid E_A=1, R_A=1) \nonumber\\
        &=& \Pr(R_A(0) = 0 \mid E_A=1, R_A(1)=1) \nonumber\\
        &=&\frac{\Pr(R_A(0) = 0, R_A(1)=1 \mid  E_A=1)}{ \Pr(R_A=1\mid  E_A=1)}\label{eq:pca}.
\end{eqnarray}

While the denominator of \eqref{pca} is estimable from data on $(E,R)$
from individuals similar to Anders, the numerator is not so easily
tamed, since, in accordance with the fundamental problem of causal
inference, we can never have data on the pair $(R(0),R(1))$.  However,
as described in \citet{apd/mm/rm}, knowledge of the joint distribution
of $(E,R)$ can yield bounds on this numerator, and thus on \pc;
furthermore, these bounds can be tightened if we can observe
additional variables, such as a covariate or a mediator.

\section{Data}
\label{sec:data}
Suppose then that we have data on $E$ and $R$ for a population of
individuals we can regard as similar to (exchangeable with) Anders, so
that probabilities estimated from the data are also relevant to him.
We initially assume that there is no confounding between $E$ and $R$,
for example because $E$ was assigned by an experimental intervention,
both for Anders and in the data set.  This is formally represented as
independence between exposure $E$ and the pair of potential responses
$\bR = (R(0),R(1))$\footnote{It is in fact enough that $E$ be
  independent of $R(0)$.}.  Then we can estimate
$${\Pr}(R(1)=1\mid  E=1) =\Pr(R=1 \mid E = 1)$$
and (using non-confounding)
$${\Pr}(R(0)=1\mid E=1)={\Pr}(R(0)=1\mid E=0) =\Pr(R=1 \mid E = 0).$$
That is, we can estimate the margins of the bivariate distribution of
$(R(0), R(1))$ given $E=1$ which appears in the numerator of
\eqref{pca}.

% In particular we can estimate
% $$\frac{\Pr(R(1) = 1 \mid E=1)}{\Pr(R(0)=1) \mid E=1)} = 
% \frac{\Pr(R = 1 \mid E = 1)}{\Pr(R=1 \mid E = 0)},$$ the {\em
%   (experimental) risk ratio\/}, $\rr$.

\section{Basic inequality}
\label{sec:basic}
The following inequality, where the
bounds are estimable, then applies \citep{apd:aberdeen}:

\begin{equation}
  \label{eq:simple-analysis}
  \max\left\{0,1 - \frac{1}{\rr}\right\} \leq  \PC_A  \leq  \min\left\{1,\frac{\Pr(R=0  \mid E = 0)}{\Pr(R=1 \mid E = 1)}\right\}.
\end{equation}
Here $\rr = {\Pr(R = 1 \mid E = 1)}/{\Pr(R=1 \mid E = 0)}$, the {\em
  (experimental) risk ratio\/}.  In particular, if $\rr>2$ the
Probability of Causation will exceed $50\%$.  In a civil court this is
often taken as the criterion to assess legal responsibility ``on the
balance of probabilities''.  {However}, if $\rr<2$ we can {not} infer
from \eqref{simple-analysis}
that $\pc_A < 50\%$.\\

\begin{example}

  \rm
  \label{ex:data}
  Suppose we have data as in \tabref{aspirin}.
  \begin{table}[htbp]
    \caption{Experimental data}
    \label{tab:aspirin}
    % \begin{center}
    \begin{tabular}[t]{lccc}
      & Myocarditis ($R=1$) & No myocarditis ($R=0$) & Total\\
      Vaccinated ($E=1$) & 30 & 70 & 100\\
      Unvaccinated ($E=0$) & 12 & 88 & 100
    \end{tabular}
      % \end{center}
  \end{table}

  To bound \pc, consider \tabref{aspirin2}.
    
  \begin{table}[htbp]
    \caption{Potential responses}
    \label{tab:aspirin2}
    % \begin{center}
    \begin{tabular}[t]{c|cc|c}
      \multicolumn{1}{c}{} &\multicolumn{2}{c}{$R(0)$}\\
      % \hline
      $R(1)$ & 0 & 1 \\
      \hline
      0 &  $88-x$ & $x-18$  & 70\\
      1 & $x$ & $30-x$  & 30\\
      \hline
                           & 88 & 12 & 100
    \end{tabular}
    % \end{center}
  \end{table}

  Since, \eg, $\Pr(R(0)=1) = \Pr(R=1 \mid E=0)$, the marginal totals
  in \tabref{aspirin2} can be copied from the internal entries of
  \tabref{aspirin}.  The internal entries of \tabref{aspirin2} can not
  be determined exactly, but have one degree of freedom, represented
  by the unspecified value of $x$.  Nevertheless, the constraint that
  every such entry must be non-negative implies that
  $18 \leq x \leq 30$.  Since (using non-confounding)
  $\pc = \Pr(R(0)=0, R(1)=1)/\Pr(R(1)=1) = x/30$, $\pc$ must lie
  between $18/30 = 60\%$ and (uninformatively) $30/30 = 100\%$.  This
  is in agreement with formula~\eqref{simple-analysis}, which indeed
  is derived in just this way.

\end{example}

\section{Complete mediator}
\label{sec:med}
%\todo{TO DO}

\begin{figure}[htbp]
 \centering
  \includegraphics[width=.4\linewidth]
    {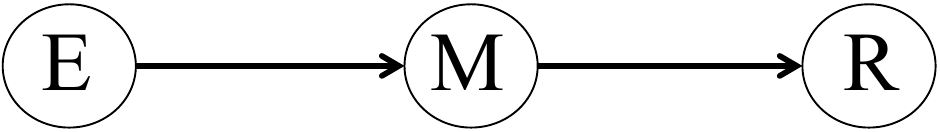}
  \caption{Complete mediator}
  \label{fig:compl-med}
\end{figure}

Consider now a situation involving a third variable $M$ that
completely mediates the causal effect of $E$ on $R$, so that $E$ has
no direct effect on $R$.  This situation is illustrated in
\figref{compl-med}.

We now introduce $M(e)$, the potential value of $M$ when $E$ takes
value $e$, and $R(m)$, the potential value of $R$ when $M$ takes value
$m$.  So the observed valued of $M$ is $M(E)$, and that of $R$ is
$R(M)$.  We assume no confounding for both the exposure-mediator and
the mediator-outcome relationship, represented formally by mutual
independence between $E$, $\bM = (M(0),M(1))$, and
$\bR = (R(0),R(1))$.  In particular this implies the Markov property
$R \cip E \mid M$, which can be tested in the data.

\cite{apd/rm/mm:sis} provide the following estimable bounds for
$\PC_A$ for the case that $M$ is observed in the experimental data (so
we know the joint distribution of $(E,M,R)$), but is not observed for
Anders:

\begin{equation} \label{eq:pc-compl-med}
\max\left\{0,{1- \frac{1}{\rr}}\right\} \leq \PC_A \leq \min\left\{0,\frac{N}{\Pr(R=1 \mid E= 1)}\right\},
\end{equation}
where the numerator $N$ in the upper bound is given by
\tabref{compl-med-upper}, with $a=\pr(M=0 \mid E=0)$,
$b=\pr(M=1 \mid E=1)$, $c=\pr(R=0 \mid M=0)$ and $d=\pr(R=1\mid M=1)$.

\begin{table}[htbp]
  \caption{Upper bound numerator, $N$}
  \label{tab:compl-med-upper}
%  \begin{center}
  \begin{tabular}{l|cc}
      & $a\leq b$ & $a>b$\\[1ex]
      \hline
      $c\leq d$ &  $a c+(1-b)(1-d)$ & $b c+(1-a)(1-d)$\\[1ex]
      % \quad\\
      $c>d$ &  $a d+(1-b)(1-c)$ & $b d+(1-a)(1-c)$
    \end{tabular}
%\end{center}
\end{table}

As we can see by comparison with \eqref{simple-analysis}, knowing
about a mediator does not improve the lower bound.  However it can be
shown that the upper bound in \eqref{pc-compl-med} is never greater
than, and is typically less than, that of (\ref{eq:simple-analysis}),
which ignores the existence of the mediator $M$.  So taking account of
a complete mediator will typically refine the bounds on \pc.\\

\begin{example}
  \rm
  \label{ex:med}
  Assume the the following probabilities:
  \begin{eqnarray*}
    \Pr(M=1 \mid E= 1)&=&0.75\\
    \Pr(M=1 \mid E = 0)&=&0.975\\
    \Pr(R=1 \mid M = 1)&=&0.1\\
    \Pr(R=1 \mid M = 0)&=&0.9
  \end{eqnarray*}
  These are consistent with the values of $\pr(R\mid E)$ given in
  \tabref{aspirin}.
  
  With these probabilities, \eqref{pc-compl-med} produces bounds
  $0.60 \leq \pc \leq 0.76$.

  If however we were to ignore the mediator $M$, we would just use
  formula \eqref{simple-analysis}, which is applicable since we still
  have no confounding between $E$ and $R$.  This yields bounds
  $0.6 \leq \pc \leq 1$, as in \exref{data}.  In such a case, ignoring
  the mediator gives a greater upper bound---less informative, though
  not inconsistent with, the more precise result obtained when taking
  it into account.
\end{example}

\section{Sufficient covariate}
\label{sec:suffcov}
\begin{figure}[htbp]
  \centering
  \includegraphics[width=.4\linewidth]
  {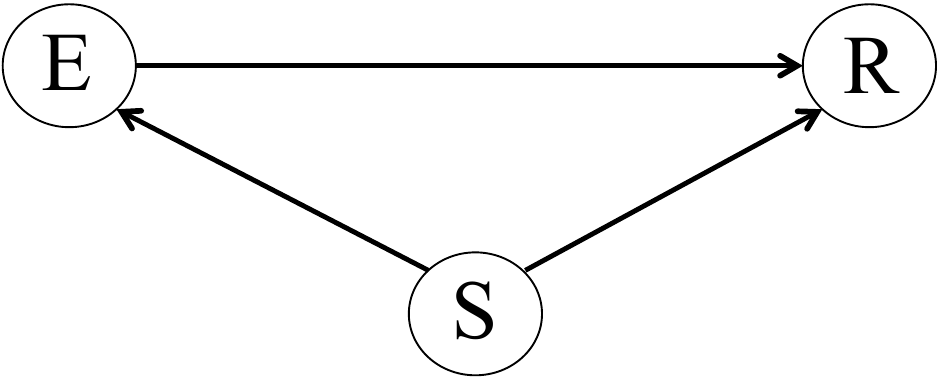}
  \caption{Sufficient covariate}
  \label{fig:suffcov}
\end{figure}

We now remove the ``no confounding'' condition, so that the above
results do not directly apply.  Instead, we suppose that there is an
additional {\em sufficient covariate\/} $S$, that can be measured
before exposure is determined and that can have an effect on both
exposure and outcome; but such that, conditional on $S$, there is no
further confounding.  This situation can be described diagrammatically
by \figref{suffcov}. 
  
If we can observe $S$ both in our data and for
Anders, the same analysis as in \secref{basic} can be applied, with
bounds as in \eqref{simple-analysis} but using probabilities further
conditioned on Anders's own value for $S$.  So now we suppose we can
observe $S$ for our data (and so estimate the joint distribution of
$(E,S,R)$), but not for Anders.

% Assume $E$ and $R$ binary and  $S$ discrete with $\Pr(S=s) >0  \quad \forall s$
 
% Suppose that from the study data we can identify 
% $$\Pr(S=s) \quad \Pr(E=e \mid S=s) \quad \Pr(R=r \mid S=s, E=e)$$

We introduce potential variables as follows.

\begin{itemize}
\item $E(s)$ the potential exposure when $S$ takes value $s$
\item $R(s, e)$ the potential response when $S$ takes value $s$ and
  $E$ takes value $e$
\end{itemize}

Given our assumption that there is no further confounding after
conditioning on $S$, we have for instance
$$ \Pr(E(s)=e) = \Pr(E=e \mid S=s).$$

In this case we have new bounds for \pc \citep{apd:aberdeen}:

  \begin{equation}
    \label{eq:better}
    \frac{\Delta}{\Pr(R= 1 \mid E=1)} \leq \pc \leq 1 -
    \frac{\Gamma}{\Pr(R = 1 \mid E = 1)}
  \end{equation}
  
  where
  \vspace{-1ex}
  \begin{align}
    \Delta = &\sum_s  \Pr(S=s | E=1)\, \times  \nonumber \\
    &{}\max\left\{0, \Pr(R= 1 \mid E = 1, S=s)
      -\Pr(R = 1 \mid E = 0, S=s)\right\} \nonumber \\
    \Gamma =  &\sum_s \Pr(S=s|E=1) \,\times  \nonumber \\
    &{}\max\left\{0, \Pr(R = 1 \mid E = 1, S=s) - \Pr(R = 0 \mid
      E= 0, S=s)\right\}. \nonumber
  \end{align}

  \subsection{Biased analysis}
  \label{sec:bias}
  Suppose we were to ignore $S$, wrongly believing that there was no
  confounding.  Then we would use the formula (2), but this would be
  incorrect.

  To see how much difference such a biased analysis can make, consider
  the following simple example, with all variables binary.

\begin{example}
  \label{ex:feels}\rm
  \hfill\\
  
  Marginally, $\pr(S=1) = \pr(S=0) = 0.5$.\\

  The distribution of $E$ given $S$ is as in \tabref{se}.
  \begin{table}[htbp]
    \centering
    \caption{$E$ given $S$}
    \begin{tabular}[t]{lccc}
      &$E=1$ & $E=0$ & Sum\\
      $S=1$ & 0.2 & 0.8 & 1\\
      $S=0$ & 0.8 & 0.2 & 1\\
    \end{tabular}
    \label{tab:se}
  \end{table}

  The distribution of $R$, given $E$ and $S$, is as in \tabref{gene}.
  \begin{table}[htbp]
    \caption{$R$ given $E$ and $S$}
    \label{tab:gene}
    % \begin{center}
    \begin{tabular}[t]{lccc}
      \multicolumn{4}{c}{$S=1$}\\
      \hline
      & $R=1$ & $R=0$ & Sum\\
      $E=1$ & 0.2 & 0.8  & 1\\
      $E=0$ & 0.8 & 0.2 & 1
    \end{tabular}\quad\quad\quad
    \begin{tabular}[t]{lccc}
      \multicolumn{4}{c}{$S=0$}\\
      \hline
      & $R=1$ & $R=0$ & Sum \\
      $E=1$ & 0.8 & 0.2 & 1 \\        
      $E=0$  & 0.2 & 0.8 & 1
    \end{tabular}
    % \end{center}
  \end{table}
  Applying the correct bounds of \eqref{better} we get
  \begin{equation}
    \label{eq:1}
    0.71 \leq \pc \leq 1.
  \end{equation}

However, suppose we ignored $S$, and incorrectly assumed no
confounding.  The distribution of $R$ given only $E$ is as in
\tabref{wrong}.
 \begin{table}[htbp]
    \centering
    \caption{$R$ given $E$}
    \begin{tabular}[t]{lccc}
      &$R=1$ & $R=0$ & Sum\\
      $E=1$ & 0.68 & 0.32 & 1\\
      $E=0$ & 0.68 & 0.32 & 1\\
    \end{tabular}
    \label{tab:wrong}
  \end{table}
  Applying the incorrect formula \eqref{simple-analysis} yields
  \begin{equation}
    \label{eq:2}
    0 \leq \pc \leq 0.47
  \end{equation}

  The contrast between the correct interval \eqref{1} and the
  incorrect interval \eqref{2} is striking.
\end{example}

\section{Mediator with sufficient covariate}\label{sec:medsuff}
\begin{figure}[h]
  \centering
  \includegraphics[width=.5\linewidth]
  {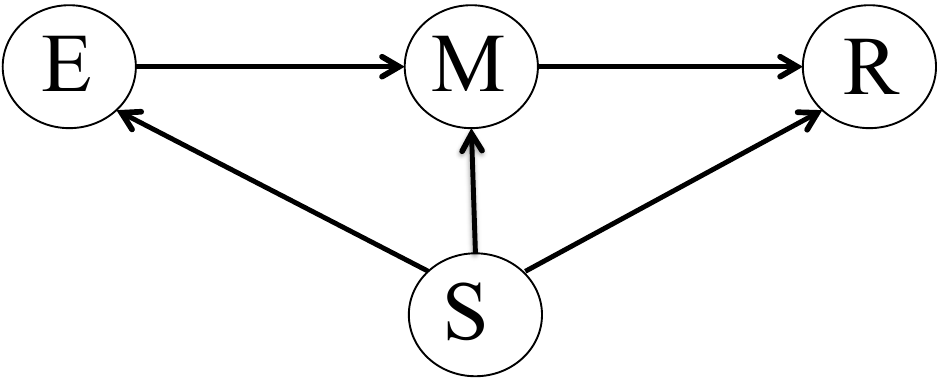}
  \caption{Complete mediator with covariate}
  \label{fig:allcomplete}
\end{figure}

\figref{allcomplete} illustrates a case involving both a mediator,
$M$, and a sufficient covariate, $S$, where $S$ can affect all of $E$,
$M$, and $R$, and conditionally on $S$ there is no direct effect of
$E$ on $R$ unmediated by $M$.  Here $E$ might denote vaccination
status; $M$, blood thickened?; $R$, thrombosis?; and $S$, sex.  As
before we assume no further confounding.

The formulas now required for the lower and upper bounds of \pc are
given in \citet[\S17.7]{apd/mm:sef}.\\

\begin{example}\rm
  \label{ex:medsuff}
  The overall distribution is defined by the following ingredients:

  \begin{eqnarray*}
    \pr(S=1) &=&	0.9\\
    \\
\pr(E=1 | S = 0)  &=&	0.9\\
    \pr(E=1 | S = 1)  &=&	0.1\\
    \\
\pr(M=1 | E=0, S=0)  &=&	0.1\\
\pr(M=1 | E=0, S=1)  &=&	0.8\\
\pr(M=1 | E=1, S=0)  &=&	0.3\\
    \pr(M=1 | E=1, S=1)  &=&	0.8\\
    \\
\pr(R=1|M=0,S=0)  &=&	0.8\\
\pr(R=1|M=0,S=1)  &=&	0.9\\
\pr(R=1|M=1,S=0)  &=&	0.7\\
\pr(R=1|M=1,S=1)  &=&	0.3
  \end{eqnarray*}

  We consider a case in which we have observed $E=R=1$, but have not
  observed $M$ or $S$.  Then using the correct formula which takes %
  into account the existence of $S$ and $M$ and the overall joint
  distribution, we obtain
  \begin{equation}
    \label{eq:3}
    0 \leq \pc \leq 0.21.
  \end{equation}
  If we ignore the existence of $M$, but still take account of $S$
  (which is not incorrect, since there is still no residual
  confounding), so conducting an analysis as in \secref{suffcov} based
  on the joint distribution of $(S,R,E)$, we get
  \begin{equation}
    \label{eq:4}
    0 \leq \pc \leq 0.53
  \end{equation}
---not wrong, but less informative than before
  because not all available information has been used.

  On the other hand, if we take account of $M$ but ignore $S$, so
  conducting a biased analysis, as in \secref{med}, based on the joint
  distribution of $(R,M,E)$, incorrectly assuming no further
  confounding and no direct effect of $R$ on $E$, we obtain
  \begin{equation}
    \label{eq:5}
    0.24 \leq \pc \leq 0.59,
  \end{equation}
totally inconsistent with the correct bounds of \eqref{3}.

Finally, on conducting a biased analysis that ignores both $M$ and
$S$, as in \secref{basic}, we get
\begin{equation}
  \label{eq:6}
   0.29  \leq   \pc  \leq  0.97,
\end{equation}
again inconsistent with \eqref{3}.
\end{example}

    \section{Conclusions}
\label{sec:conc}

Even when presented with substantial unconfounded experimental data
fully determining the probabilistic dependence of a response variable
on an exposure variable, we are typically only able to establish
bounds for the probability of causation in a case where an individual
has developed the response after being exposed.  However these bounds
can be enhanced or adjusted in the presence of additional information,
such as data on covariate or mediator variables.  Incorporating
mediator variables into the analysis will typically narrow the bounds
for the probability of causation, by shedding light on the mechanisms
through which the exposure affects the response.  Taking mediator
variables into account can improve the accuracy of causal inference,
but ignoring them still yields valid conclusions, though these will be
less precise.

In general, however, the data employed in cause-effect problems are
observational, and susceptible to various forms of distortion: for
instance, there may be confounding variables, which it is crucial to
consider in order to ensure the reliability of our causal conclusions.
We have shown, by examples, that when confounding variables are present
but not observed, or ignored, an analysis that wrongly assumes no
confounding can lead to biased inferences, totally inconsistent with
the correct conclusions.

Further work is required to investigate more complex scenarios and
real-world applications.

%\bibliography{strings,causal}
\end{document}